\documentclass[11pt,a4paper]{amsart}
\usepackage{amscd,amssymb,amsmath,latexsym}

\theoremstyle{plain}
\newtheorem{theorem}{Theorem}[section]
\newtheorem{lemma}[theorem]{Lemma}
\newtheorem{proposition}[theorem]{Proposition}
\newtheorem{corollary}[theorem]{Corollary}

\theoremstyle{definition}

\newtheorem{example}[theorem]{Example}

\theoremstyle{remark}

\newcommand{\bfx}{{\boldsymbol{x}}}
\newcommand{\bfu}{{\boldsymbol{u}}}
\newcommand{\bfU}{{\boldsymbol{U}}}
\newcommand{\bfy}{{\boldsymbol{y}}}
\newcommand{\bfzeta}{{\boldsymbol{\zeta}}}
\newcommand{\bfxi}{{\boldsymbol{\xi}}}

\renewcommand{\epsilon}{\ensuremath{\varepsilon}}
\newcommand{\C}{\ensuremath{\mathbb{C}}}

\newcommand{\Q}{\ensuremath{\mathbb{Q}}}
\newcommand{\N}{\ensuremath{\mathbb{N}}}

\newcommand{\Z}{\ensuremath{\mathbb{Z}}}

\newcommand{\h}{{\rm h}}
\newcommand{\Ch}{{\rm Ch}}

\makeatletter
\newif\if@borderstar
\def\bordermatrix{\@ifnextchar*{%
  \@borderstartrue\@bordermatrix@i}{\@borderstarfalse\@bordermatrix@i*}%
}
\def\@bordermatrix@i*{\@ifnextchar[{%
  \@bordermatrix@ii}{\@bordermatrix@ii[()]}
}
\def\@bordermatrix@ii[#1]#2{%
  \begingroup
    \m@th\@tempdima8.75\p@\setbox\z@\vbox{%
      \def\cr{\crcr\noalign{\kern 2\p@\global\let\cr\endline }}%
      \ialign {$##$\hfil\kern 2\p@\kern\@tempdima & \thinspace %
      \hfil $##$\hfil && \quad\hfil $##$\hfil\crcr\omit\strut %
      \hfil\crcr\noalign{\kern -\baselineskip}#2\crcr\omit %
      \strut\cr}}%
    \setbox\tw@\vbox{\unvcopy\z@\global\setbox\@ne\lastbox}%
    \setbox\tw@\hbox{\unhbox\@ne\unskip\global\setbox\@ne\lastbox}%
    \setbox\tw@\hbox{%
      $\kern\wd\@ne\kern -\@tempdima\left\@firstoftwo#1%
        \if@borderstar\kern2pt\else\kern -\wd\@ne\fi%
      \global\setbox\@ne\vbox{\box\@ne\if@borderstar\else\kern 2\p@\fi}%
      \vcenter{\if@borderstar\else\kern -\ht\@ne\fi%
        \unvbox\z@\kern-\if@borderstar2\fi\baselineskip}%
        \if@borderstar\kern-2\@tempdima\kern2\p@\else\,\fi\right\@secondoftwo#1 $%
    }\null \;\vbox{\kern\ht\@ne\box\tw@}%
  \endgroup
}
\makeatother

\def\Qrng{\Q[\bfx]}
\def\Zrng{\Z[\bfx]}
\def\Crng{\C[\bfx]}
\def\mod{\,\textup{mod}\,}

\title[Effective bounds for polynomial systems]{Effective bounds for polynomial systems \\defined over
the rationals}
\author[T. Krick]{Teresa Krick}
\address{Departamento de Matem\'atica, FCEN, Universidad de Buenos Aires and IMAS, UBA-CONICET,  Argentina.}
\email{krick@dm.uba.ar}
\urladdr{http://mate.dm.uba.ar/\~\,krick}

\begin{document}

\begin{abstract}  This survey article aims to present some material that illustrates the kind of estimates one can obtain in effective algebraic geometry, for affine polynomial equation systems defined over the rational numbers, and focuses on the case of finite varieties.
\end{abstract}


\maketitle

\vskip-0.5cm
\begin{center}{\em Dedicated to my advisor, Joos Heintz, 27 October 1945--3 October 2024}\end{center}

\section{Introduction}

These pages aim to present some material that illustrates the kind of estimates we can obtain in  effective algebraic geometry, for affine polynomial equation systems defined over the rational numbers $\Q$. It is very much inspired by the work and interests of Joos Heintz, contained in his PhD thesis and in later work and methods he promoted in the  nineties with the members of his working group Noa\"i  Fitchas I belonged to.

 I'll consider the case where polynomials are assumed to be given in {\em dense representation}: a polynomial $f$ in $n$ variables is described by an a priori bound $d$ for its {\em degree} and the array of its $\binom{d+n}{n}\le (n+1)^d$ coefficients (zero-coefficients as well as non-zero ones), and when it has integer coefficients, a bound $h$ for their {\em heights}, i.e.  their bit sizes. The accent here is put on describing  degree and height bounds for polynomials and polynomial identities  arising when solving some classical  problems, not on the algorithms to obtain them. Bounds are interesting per se as a complexity measure, but  are also  useful to know in advance when performing any algorithm, or   to certify results obtained by numerical algorithms. Although   the estimates look somehow coarse, I preferred  to present  them in detail instead of adopting the ${\mathcal O}$ notation, to show that the bounds are precise with definite universal constants.

Wishing to stay at an understandable survey level,  I'll only focus here on systems of equations defining a finite (possibly empty) number of solutions in  $\C^n$: A  friendly introductory reference that presents this setting is the famous book by Cox, Little and O'Shea \cite[Ch.4 \& 5]{CLO2015}. Needless to say,
there are in the literature many more results and  extensions of what I am presenting here, obtained by a huge community, in multiple different settings and directions:
  \begin{itemize} \item Algorithms,
  \item for positive dimension systems of equations as well,
  \item  in the projective, multiprojective,  toric settings,
  \item with their corresponding encodings and measures of the data, including  the {\em straight-line program} and {\em sparse} encodings.
  \end{itemize}
This is just a first step to give some flavor of  problems that one can consider, and what is done or can be done in effective algebraic geometry over the rationals.  I apologize for the biased choice of references and to the people who have done work that I do not cite.\\

{\em Acknowledgmwents.}
The material included here was presented as a course at the Journ\'ees Nationales de Calcul Formel 2025. I want to express my deep gratitude  to the organizers  for inviting me to deliver it and  for giving me the opportunity to reexplore it. I am  also especially indebted to Mart\'\i n Sombra for his very careful reading and detailed suggestions, and to Nicol\'as Allo G\'omez, Alin Bostan, Pierre Lairez, Gr\'egoire Lecerf  and the referees for their interest and comments, which  all contributed to improve the presentation.  Santiago Laplagne, Mart\'\i n Mereb and Tushant Mittal  helped me figure out a proof for Example~\ref{ex:perron}.  Last but not least, several results described here are slight variants of developments in recent joint work  \cite{BKM2024} with  Lorenzo Baldi and Bernard Mourrain: many thanks to them too!

\section{Setting and notations}

In all what follows, the set of natural numbers is for me, as in Argentina and several  other places in the world, the set $\N=\{1,2,3,\dots\}$  (without $0$)   and $\N_0:=\N\cup \{0\}$, $\Z $ is the ring of integer numbers, and $\Q$  and $\C$ are the fields of rational  and complex numbers respectively.\\ Given $n,s\in \N$,  $\bfx=\{x_1,\dots,x_n\}$ is a set of $n$ variables, and
$$f_1,\dots, f_s\in \Crng=\C[x_1,\dots,x_n]$$ are $s$ (nonzero) polynomials such that  the (affine) algebraic variety
\begin{equation}\label{eq:v}V:=V_\C(f_1,\dots,f_s)=\{\bfzeta \in \C^n: f_1(\bfzeta )=\dots =f_s(\bfzeta )=0\}\ \subset \, \C^n\end{equation}
 is finite: either empty or {finite} of cardinality $\ge 1$. \\
  In terms of the ideal $I=(f_1,\dots,f_s)\subset \Crng$,  by the Nullstellensatz (see Section \ref{sec:NSS}), $$V=\emptyset \ \iff \ I=\Crng,$$ while when $V$ is finite and non-empty,  we say that $I$ is a   {\em zero-dimensional ideal} (or $V$ is  a zero-dimensional  variety), which is equivalent to the
fact that   the $\C$-vector space $\Crng/I$ is of finite dimension $\ge 1$. In this last case,
we set $\deg(V):=\# (V)$ and $D:=\dim_\C\big(\Crng/I\big)$.

\smallskip

For zero-dimensional ideals, we necessarily have, by Krull's height theorem, that $s\ge n$ since $k$ polynomials in $\Crng$ define a variety, which --if not empty-- is  of dimension at least $n-k$.

\smallskip

We will soon consider the case when the  variety $V\subset \C[\bfx]$ introduced in \eqref{eq:v} is defined over the rational numbers $\Q$, i.e. we take  $f_1,\dots,f_s\in \Qrng$, or even simpler,  in $\Zrng$, that is, with integer coefficients, and $I=(f_1\dots,f_s)\subset \Qrng$.  We observe that the dimension $D$ remains the same when considering the quotient rings $\Qrng/I$ or $\Crng/I$, where $I\subset \Qrng$ or $\Crng$ accordingly.

\smallskip
The (logarithmic) {\em height} $\h(f)$ of any non-zero polynomial $f$ with integer coefficients, in any number of variables, is here defined --as usual in computational algebra--  as the maximum of the binary  logarithms of the absolute values of its (nonzero) coefficients, which are all integers: If
$$f=\sum_{\alpha \in \N_0^n} c_\alpha \bfx^\alpha\in \Zrng, $$
then
$$ \h(f):=\max_\alpha\{\log | c_\alpha| : c_\alpha\ne 0\},$$
with the convention that in this text, $\log=\log_2$. \\ This height essentially coincides with the maximum bit size of all the coefficients and is another measure of the input polynomials that we will use in this text, in addition to the degrees of these polynomials, their number of variables, and if needed the number of polynomials itself.

\section{The B\'ezout inequality}
The B\'ezout inequality we are going to use all along is an affine  version of the  B\'ezout theorem in projective varieties over algebraically closed fields.

\medskip
\begin{theorem}{\bf (B\'ezout inequality)}\label{thm:BI} \\
Let $f_1,\dots,f_s\in \Crng\setminus \{0\}$ define a zero-dimensional variety $V\subset \C^n$.\\ Set  $d_j:=\deg(f_j)$ and assume that   $d:=d_1\ge d_2\ge \cdots \ge d_{s-1}$ (with no condition with respect to $d_s$). Then
$$ \deg( V)\le d_1\cdots d_{n-1} d_s\le d^{n-1}d_s.$$
\end{theorem}

\smallskip
(In this theorem it clearly makes sense to take $f_s$ as the polynomial of smaller degree, but  I stated it this way for coherence with Theorem~\ref{prop:ABI} below.)

\bigskip

For example, taking $d,h\in \N$, for $$V_1=V_\C(x_1^d-1, \dots, x_n^d-1) \ \subset \C^n,$$
 we have $V_1=\{\bfzeta=(\zeta_1 , \dots, \zeta_n)\in \C^n: \, \zeta_1^d=\dots =\zeta_n^d=1\}$ with $\deg( V_1)=d^n$, while for
 \begin{equation}\label{eq:V2}V_2=V_\C(x_1-2^h,x_2-x_1^d, \dots , x_n-x_{n-1}^d)\ \subset \C^n,\end{equation}  we have $V_2=\{(2^h,2^{dh},\dots, 2^{d^{n-1}h})\}$ with $\deg( V_2)=1<  d^{n-1}$ when $d>1$.

\medskip
This  seemingly lighter version of the classical B\'ezout theorem is not a special case of it, since there are zero-dimensional systems that have a finite number of common zeros in $\C^n$ but an infinite number when taking their zeros at infinity; consider for instance the polynomials $$f_1=xy, f_2=xz, f_3=x-1\ \in \C[x,y,z]$$ which satisfy that $V=\{ (1,0,0)\}$, although if one homogenizes  the polynomials with a first variable, then there are infinite projective zeros $$(1:1:0:0)\cup (0:0:0:1) \cup \{(0:0:1:z)\,:\, z\in \C\}.$$

There are several proofs for this B\'ezout inequality in the literature, see for instance \cite[Ch.2, Cor.]{MaWh1983}, which also counts  multiplicities, or \cite[Ex.8.4.6]{Ful1984}, as a consequence of an extension of the classical B\'ezout theorem to arbitrary projective varieties. I will sketch here how it can be
easily deduced from  Joos Heintz' PhD thesis, published in \cite{Hei1983}: There, it is shown in Thm.1 that if $V,W\subset \C^n$ are arbitrary affine algebraic varieties, then $$\deg(V\cap W)\le \deg(V)\deg(W),$$ where $\deg(V)$ is the (geometric) {\em degree} of the variety $V$, which for an irreducible variety of dimension $r$ is the maximum finite number of  points that one can obtain when intersecting the variety with $n-r$ affine hyperplanes, and  for an arbitrary variety  is defined as the sum of the degrees of its irreducible components; for a hypersurface $V=V(g)$ defined by a squarefree polynomial  $g$, we have $\deg(V)=\deg(g)$, and for a zero-dimensional variety $V$, $\deg(V)=\# (V)$.
\\
To deduce that in the B\'ezout bound, the degrees of  only the  $n$ polynomials $f_1,\dots,f_{n-1},f_s$ suffice,  one can  show that there exist ``upper triangular" linear combinations of the input polynomials $f_1,\dots, f_s$ of the form
$$\left\{\begin{array}{ccccccccccccc}
g_1& =&f_s&  {\ }&  &  &  &   &   &   & & &\\
g_2& =& & &a_{2,1}f_1 &+ &\cdots & + & a_{2,n-1}f_{n-1} & + &\cdots &+&a_{2,s-1}f_{s-1}\\
\vdots&  &&&&&\ddots &&&&&&\\
 g_n&=&&&&&&& a_{n,n-1}f_{n-1}& + &\cdots& +& a_{n,s-1}f_{s-1},\end{array}\right.$$
 for $a_{i,j}\in \C$,
such that $\dim\big(V_\C(g_1,\dots,g_k)\big)=n-k$ and therefore, $V_\C(g_1,\dots,g_n) $,  which obviously contains $V$, is finite (see for instance \cite[Proof of Thm.14]{CGH1989}).\\
Note that it is not true that there always exist  $n$ such linear combinations $g_1,\dots, g_n$ so that $V_\C(g_1,\dots,g_n)=V$, as shown by the very simple example
$(x^3-x, x^2-2x)\subset \C[x]$.

\smallskip
The B\'ezout inequality proves to be very useful in Computer Algebra, because it doesn't make any assumption on the polynomials, neither the way they intersect (no consideration of {\em genericity} is needed) nor their number.

\section{An Arithmetic B\'ezout inequality}\label{sec:2}

The {\em (logarithmic) height} $\h(V)$ of a zero-dimensional  variety $V$ defined over $\Q $ is a parameter that essentially measures the complexity  of the points in $V$, and in particular yields a bound for the absolute values of their coordinates. It is a special case of a more general theory of heights of varieties of any dimension. I won't use the most sophisticated definition of height, but rather a binary version of the notion introduced by Patrice Philippon in \cite{Phi1991,Phi1995}, which is sufficient for our purposes, and is defined through the {\em Chow form} of $V$. A general presentation of the Chow form (for varieties of any dimension over any field) can be found in \cite[Sec.I.6.5]{Sha1974}, and for the connection with the height of a variety I follow  \cite[Sec.1]{KPS2001}.

\smallskip

Since $V$ is defined over $\Q$, if we  denote $\bfzeta=(\zeta_1,\dots,\zeta_n)$ for $\bfzeta\in V\subset \C^n$,  the polynomial
$$ \prod_{\bfzeta\in V}(U_0+\zeta_1U_1+\dots+\zeta_nU_n)$$ in $n+1$ new variables $\bfU=(U_0,\dots U_n)$ happens to belong to $\Q[\bfU]$ and we can choose $c\in \N$  such that
\begin{equation}\label{eq:chow}\Ch_V(\bfU):=c\prod_{\zeta\in V}(U_0+\zeta_1U_1+\dots+\zeta_nU_n)\ \in \,\Z[\bfU]\end{equation}
 is a primitive polynomial in $\Z[\bfU]$: this  is the  (primitive) Chow form of the zero-dimensional variety $V$. It is an homogeneous polynomial of degree $\deg( V)$, which satisfies that for any nonzero $\bfu\in \C^{n+1}$,
\begin{align*}\Ch_V(\bfu)=0 & \iff u_0+\zeta_1u_1+\dots +\zeta_nu_n=0 \mbox{ for some }  \bfzeta\in V \\
& \iff \exists\, \bfzeta \in V :\ \bfzeta  \mbox{ belongs to the hyperplane } \{u_1x_1+\dots + u_nx_n=-u_0\}.\end{align*}
The Chow form of a variety (of any dimension) characterizes the variety in the sense that different varieties have different Chow forms, and one can recover the variety from its Chow form.\\

The {\em (Philippon-)height} $\overline\h(V)$ of our zero-dimensional variety $V$ defined over $\Q$ is  $$\overline\h(V):=\int_{S_{n+1}}
\ln |\Ch_V| \ \mu_{n+1}+\big(\sum_{i=1}^n\frac{1}{2i}\big) \deg(V) $$ (here  $S_{n+1}:=\{ (z_0,\dots,z_n)\in
\C^{n+1}: |z_0|^2+\cdots+|z_n|^2=1\}$ is the unit sphere in
$\C^{n+1}$, $\ln$ is the natural logarithm   and $\mu_{n+1}$ is the measure on $ S_{n+1}$ of total mass $1$, invariant
with respect to the unitary group $U(n+1)$). Again this is a special case of a definition of  height for varieties of any dimension defined over any number field (see for example \cite{Phi1995} or \cite{KPS2001}).  Since we are interested here in bit sizes (binary logarithms) rather than natural logarithms, I  use here a slight variant of $\overline\h(V)$, namely $\h(V):=\log e   \overline \h(V)$. In our case  of a zero-dimensional variety $V$ it  satisfies the following properties
\begin{align}
&  \sum_{\bfzeta \in V} \log(\|(1,\bfzeta)\|_2)\le \h(V) , \  \mbox{ where } \| \cdot\|_2 \ \mbox{denotes the Euclidean norm}; \label{eq:rootIn} \\
 &
|\h(V)-\h(\Ch_V)|\le 3 \log(n+1)\deg (V).  \label{eq:hVChV}\end{align}
(For \eqref{eq:rootIn} see
\cite[Prop.4]{Phi1991} while \eqref{eq:hVChV} follows  combining inequalities by different  authors, as mentioned for instance in
\cite[Ineq.(1.1)\&(1.2)]{KPS2001}.)

\smallskip
I illustrate the notion of height of a zero-dimensional variety and the previous inequalities  by the toy case where  $V$ is composed by  a single rational point
$$\bfzeta =(\frac{a_1}{b_1},\dots,\frac{a_n}{b_n}) =\frac{(c_1,\dots,c_n)}{c}\in \Q^n$$ with $a_i,b_i,c_i\in \Z$ and $c\in  \N$  a minimal common denominator: In this toy example,
$$\h(V)=\log(\|(1,\frac{a_1}{b_1},\dots,\frac{a_n}{b_n} )\|_2)+\log c =\frac12 \log(c^2+c_1^2+\dots+c_n^2),$$
(see \cite[Sections 1.2.3 \& 1.2.4]{KPS2001}) while $$\h(\Ch_V)=\h(c\,U_0+c_1U_1+\cdots+c_nU_n)=\max_i\{\log c ,\log |c_i| \}.$$
For instance, for $\zeta=\frac{c+1}{c}$ with $c\in \N$, which as a rational number satisfies $\h(\zeta)=\log (c+1)$ (this somehow represents its ``complexity"), we have $ \frac12+\log c <\h(V)<\frac12+\log(c+1)$ while $\frac12\le \log(\|(1,\zeta)\|_2)\le \frac{\log 5}{2}$: This shows that Inequality~\eqref{eq:rootIn}  may not be sharp and that
 $\h(V)$ is a more accurate measure of the complexity of the point  in $V$ than just its  real size.\\

I describe now an arithmetic B\'ezout inequality in terms of the degrees and heights of the integer polynomials that define our zero-dimensional variety $V$,
which is a consequence of Corollary 2.11 in \cite{KPS2001}. Its special feature is that it is able to   distinguish the degree and height of one of the polynomials with respect to the others.

\begin{theorem} \label{prop:ABI} {\bf (An Arithmetic B\'ezout inequality)} \\Let $f_1,\dots,f_s\in \Zrng\setminus \{0\}$ define a zero-dimensional variety $V\subset \C^n$.\\ Set  $d_j:=\deg(f_j)$, $h_j:=\h(f_j)$ for $1\le j\le s$. Assume that   $d:=d_1\ge d_2\ge \cdots \ge d_{s-1}$ (with no condition with respect to $d_s$), and set $h:=\max\{h_1,\dots, h_{s-1}\}$.  Then
\begin{align*} \h( V)& \le \Big(\frac{h_s}{d_s}+ \sum_{j=1}^{n-1}\frac{h}{d_{j} }+2n\log(n+1)\Big)d_1\cdots d_{n-1}d_s \\
& \le   d^{n-1}h_s+ (n-1)d^{n-2} d_s h+2n\log(n+1)d^{n-1}d_s.
\end{align*}
\end{theorem}

\begin{proof} I give the proof because this statement doesn't appear  anywhere in this form. It makes use of the notion of height of varieties of any dimension.\\
We apply \cite[Cor.2.11]{KPS2001} to the equidimensional variety $V_\C(f_s)$ of dimension $n-1$, which satisfies
$$\deg(V_\C(f_s))\le d_s \quad \mbox{and} \quad \h(V_\C(f_s))\le h_s+(n+1)\log(n+1)d_s,$$
by application to $V_\C(f_s)$ of the comment after Corollary 2.11 in the same paper. \\
Then
\begin{align*}
\h\big(V_\C(f_1,\dots,f_s)\big) & = \h\big(V_\C(f_s)\cap V_\C(f_1,\dots,f_{s-1})\big)\\
&\le \Big( \h(V_\C(f_s)) + \big(\sum_{j=1}^{n-1} \frac{1}{d_j}\big) h\deg(V_\C(f_s)) \\ & \qquad \qquad\qquad \qquad + (n-1)\log(n+1)\deg(V_\C(f_s))\Big)\prod_{i=1}^{n-1} d_i\\
& \le \Big(  h_s+(n+1)\log(n+1)d_s \\ & \qquad \qquad\qquad \qquad+ \big(\sum_{j=1}^{n-1} \frac{1}{d_j}\big) d_sh + (n-1)\log(n+1)d_s\Big)\prod_{i=1}^{n-1} d_i.
\end{align*}
The second statement  follows directly from this, since $\big(\displaystyle\sum_{j=1}^{n-1} \frac{1}{d_j}\Big)\prod_{i=1}^{n-1} d_i\le (n-1)d^{n-2}$.
\end{proof}

\bigskip
Inequality \eqref{eq:rootIn} above and the Arithmetic B\'ezout inequality directly imply the following upper bound for the coordinates of any root $\bfzeta=(\zeta_1,\dots,\zeta_n)\in V$:
For $d_s\le d$ and $h_s\le h$ we get
\begin{equation}\label{eq:hroot}\log |\zeta_i| \le  nd^{n-1}h+ 2n\log(n+1)d^{n},\quad 1\le i\le n.\end{equation}
 It also shows that the arithmetic B\'ezout bound is quite tight in terms of the height dependence since for instance, for the variety $V_2=\{\bfzeta=(2^h,2^{dh},\dots, 2^{d^{n-1}h})\}$ defined in \eqref{eq:V2} above, we have  $$\log(\zeta_n)= d^{n-1}h\le \h(V_2)\le d^{n-1}h+\frac12 \log(n+1)$$ while the  estimate in Theorem~\ref{prop:ABI} gives
$$\h(V_2)\le d^{n-1}h +2n\log(n+1)d^{n-1}.$$

Bound \eqref{eq:hroot} for the coordinates implies in turn an upper bound for  the value of any complex polynomial $p\in \Crng$ at a root of our zero-dimensional variety $V$ defined over $\Q$:
\begin{corollary} \label{cor:hproot} {\bf (Upper bounds for the roots)}\\ Let $f_1,\dots,f_s\in \Zrng \setminus \{0\}$ define a zero-dimensional variety $V\subset \C^n$, with $\deg(f_j)\le d$ and $\h(f_j)\le h$ for $1\le j\le s$.
Let $p\in \Crng\setminus \{0\}$ and set $d_p:=\deg(p)$ and  $h_p$ for the maximum of the binary logarithms of the absolute values of the (nonzero) coefficients of $p$. Then, for every $\bfzeta \in V$ we have
$$\log |p(\bfzeta)| \le   h_p+ nd^{n-1}d_p h+ 3n\log(n+1)d^{n}d_p.$$
\end{corollary}
\begin{proof} Let $p=\sum_{|\alpha|\le d_p}a_\alpha \bfx^\alpha$, which has  $\binom{d_p+n}{n}\le (d_p+1)^{n}$ coefficients. Then, by the bound for the coordinates \eqref{eq:hroot} we have
$$\log(|\zeta_1|^{\alpha_1}\cdots |\zeta_n|^{\alpha_n})\le  nd^{n-1}d_ph+2n\log(n+1)d^{n}d_p,$$ and therefore, considering the multiplication by $|a_\alpha|$ and adding up the $(d_p+1)^{n}$ terms, we obtain
\begin{align*}\log|p(\bfzeta)|&\le  h_p+ nd^{n-1}d_p h+ 2n\log(n+1)d^{n}d_p+  n\log(d_p +1)\\ &\le  h_p+ nd^{n-1}d_p h+ 3n\log(n+1)d^{n}d_p.\end{align*}
\end{proof}

 But more interestingly, we also obtain  lower bounds for the non-zero coordinates of the roots in $V$, and also for the separation of these roots, which extend the bounds for univariate integer polynomials that appear for example in the beautiful book by Maurice Mignotte \cite{Mig1992}, and agree with those that can be found in the literature (\cite{Can1987}, \cite[Thm.45, Cor.49]{Yap2000}, \cite[Cor.9]{EMT2020}), with the advantage that we don't need to assume that we have exactly $n$ equations that determine the zero-dimensional variety, or that the system is still zero-dimensional if one considers zeros at infinity. The following is \cite[Lem.4.2 \& Lem.4.10]{BKM2024}, where we can drop the assumption that $\h(p)\le h$ or $\deg(p)\ge d$ due to the application of the  more precise Theorem~\ref{prop:ABI}. (The upper bounds are not as sharp as in the previous corollary, but I include them anyway for symmetry.)

\begin{proposition} \label{prop:LSbounds} {\bf (Lower and separation bounds for the roots)}\\Let $f_1,\dots,f_s\in \Zrng\setminus \{0\}$ define a zero-dimensional variety $V\subset \C^n$, with $\deg(f_j)\le d$ and $\h(f_j)\le h$ for $1\le j\le s$.
Let $p\in \Zrng\setminus \{0\}$ and set $d_p:=\deg(p)$ and $h_p:=\h(p)$.
 Then, for every $\bfzeta,\bfxi \in V$ we have
\begin{enumerate}\item $\big| \log  |p(\bfzeta)| \big| \le d^n h_p + 2nd^{n-1}d_ph + 4(n+1)\log(n+2)d^{n}d_p $ \  if   $p(\bfzeta)\ne 0$;\vskip2pt
\item $ \big|\log |p(\bfzeta)-p(\bfxi)| \big|\le d^{2n}h_p + 4nd^{2n-1}d_ph + 4(2n+1)\log(2n+2)d^{2n}d_p  $ \  if   $p(\bfzeta)\ne p(\bfxi)$.
\end{enumerate}
In particular,  for every $\bfzeta=(\zeta_1,\dots,\zeta_n), \bfxi=(\xi_1,\dots,\xi_n)\in V$, we have
\begin{enumerate}\item $\big| \log  |\zeta_i| \big| \le 2nd^{n-1}h + 4(n+1)\log(n+2)d^{n}$ \ if $\zeta_i\ne 0$; \vskip2pt
 \item $\big|\log  |\zeta_i-\xi_i|  \big|\le 4nd^{2n-1}h + 4(2n+1)\log(2n+2)d^{2n}$ \ if $\zeta_i\ne \xi_i$.
 \end{enumerate}
\end{proposition}
\begin{proof} To prove Item (1), we apply the arithmetic B\'ezout inequality together with Inequality \eqref{eq:rootIn} above to the zero-dimensional ideals in $n+1$ variables $$J_1=\big(f_1(\bfx),\dots,f_s(\bfx),x_{n+1}-p(\bfx)\big)\quad \mbox{and}\quad J_2=\big(f_1(\bfx),\dots,f_s(\bfx), 1-x_{n+1}p(\bfx)\big).$$
Since $V_\C(J_1)=\{(\bfzeta,p(\bfzeta)):\,\bfzeta\in V\}$ and $V_\C(J_2)=\{ (\bfzeta,1/p(\bfzeta)): \, \bfzeta\in V,\, p(\bfzeta)\ne 0\} $ we deduce
$$\log |p(\bfzeta)| \le \sum_{\bfzeta\in V}\log(\|(1,\bfzeta,p(\bfzeta))\|_2)\le d^n h_p + nd^{n-1}d_ph+ 2(n+1)\log(n+2)d^{n}d_p.$$
Similarly, \begin{align*}\log \left|\frac{1}{p(\bfzeta)}\right| &\le \sum_{\bfzeta\in V}\log(\|(1,\bfzeta,\frac{1}{p(\bf\zeta)})\|_2)\\ & \le d^n h_p+ nd^{n-1}(d_p+1)h+ 2(n+1)\log(n+2)d^{n}(d_p+1)\end{align*}
implies that $\log |p(\bfzeta)| \ge - \big(d^n h_p + 2nd^{n-1}d_ph + 4(n+1)\log(n+2)d^{n}d_p\big)$.\\
To prove Item (2), we apply Item (1) to  the zero-dimensional ideal in $2n$ variables $$J=(f_1(\bfx),\dots,f_s(\bfx),f_1(\bfy),\dots,f_s(\bfy))\subset \C[\bfx,\bfy]$$ and the polynomial $p(\bfx)-p(\bfy)$ of degree $d_p$ and height $h_p$.\\
We then obtain the bounds for the coordinates and the separation of the roots applying this to the polynomials $p=x_i$ and $p=x_i-y_i$ respectively, of degree 1 and height 0.
\end{proof}

These general lower bound for the coordinates of the roots are again essentially  tight as shows the following example presented in \cite{Can1987}:
$$V_3=V_\C(2^hx_1-1,x_2-x_1^d,\dots,x_n-x_{n-1}^d)$$
where the unique root $\bfzeta=(1/2^h,\dots, 1/2^{d^{n-1}h})$
satisfies $\log \zeta_n =-d^{n-1}h$ while the lower bound in Proposition~\ref{prop:LSbounds} gives
$$\log \zeta_n \ge - \big( 2nd^{n-1}h + 4(n+1)\log(n+2)d^{n}  \big).$$
They are also quite tight for the separation of the roots, according to  the multivariate Mignotte-type example
 developed in \cite[Sec.4]{EMT2020}.

\bigskip

Finally, applying more carefully Inequality \eqref{eq:rootIn}  to subsets of roots, these statements can be adapted to obtain more general estimates for sums of upper, lower and separations bounds for the roots, as used for example in \cite[Lem.4.10]{BKM2024}.

\section{An arithmetic Nullstellensatz} \label{sec:NSS}
The Nullstellensatz (NSS) is a cornerstone in algebraic geometry, which first appeared in a complete form in Hilbert's work \cite{Hil1893}.  These are its two  forms, known as  {\em weak NSS} and  {\em strong NSS}, which are in fact equivalent:

\begin{theorem} {\bf(The Nullstellensatz)}\\
Let $I\subset \Crng$ be an ideal. Then,
\begin{align*} (1)&\quad \mbox{Weak NSS:} \quad
V_\C(I)=\emptyset \ \iff \ 1\in I;\\
(2)& \quad \mbox{Strong NSS:} \quad I(V_\C(I))=\sqrt I.\end{align*}
\end{theorem}
Here $\sqrt I :=\{ f\in \Crng:\, \exists \, N\in \N \mbox{ s.t. } f^N\in I\}$ is the {\em radical} of the ideal $I$, and given a set $X\subset \C^n$, $I(X):=\{f\in \Crng:\, \forall \, \bfzeta \in X, \, f(\bfzeta)=0 \}$ is the {\em vanishing ideal} of the set $X$.\\

Both statements are quite well known in the univariate case but do not readily generalize to the multivariate case since  $\C[\bfx]$ is not a principal ideal domain when the number of variables   is y greater than $1$. For both statements,  the $\Leftarrow$ and $\supseteq $ directions are obvious, and (2) also easily implies (1). The passage of (1) to (2) is usually done trough  the famous {\em Rabinowitsch trick}.
There are several proofs for the (weak) NSS: one particularly elementary, based on an extension theorem and resultants, can be found in \cite[Ch.4, Thm.2]{CLO2015} for instance. Another more classical proof  using {\em Zariski's lemma} shows an equivalent statement:    Maximal ideals $\mathcal{M}\subset \Crng$ are all of the form $\mathcal{M}=(x_1-\zeta_1,\dots,x_n-\zeta_n)$ for  some $\bfzeta=(\zeta_1,\dots,\zeta_n)\in \C^n$.\\

When the ideal $I=(f_1,\dots,f_s)$ is generated by polynomials $f_1,\dots,f_s\in \Zrng$, then the weak NSS reads as
\begin{quote}{\em  Let $f_1,\dots,f_s\in \Zrng$ be polynomials such that the equation system
$$f_1(\bfx)=0,\dots,f_s(\bfx)=0$$
has no solution in $\C^n$. Then there exists $a\in \N$ and $g_1,\dots,g_s\in \Zrng$ satisfying the B\'ezout identity
$$a = g_1f_1+\cdots + g_sf_s.$$}
\end{quote}

\bigskip
As it is, this is a noneffective statement. Effective versions of the NSS estimate the degrees and the heights of polynomials satisfying the B\'ezout identity, and apply to many situations in number theory, theoretical computer science and computer algebra. I mention here the following particular case of the best current height estimate,  from \cite[Cor.4.38]{DKS2013}, which ``arithmeticizes" the (best) proof for degrees over varieties by Jelonek in \cite[Cor.1.1]{Jel2005}, after pioneering results by Hermann (1926), Seidenberg (1974), Masser and W\"ustholz (1983), Brownawell (1987), Caniglia-Galligo-Heintz (1988) and Koll\'ar (1988) for the degrees, and by Berenstein-Yger (1991) for height estimates (see also \cite{KrPa1996} and \cite{KPS2001}).

\begin{theorem}\label{thm:ANSS} {\bf (An Arithmetic Nullstellensatz)}\\
Let $f_1,\dots,f_s\in \Zrng\setminus \{0\}$ be s.t. $V_\C(f_1,,\dots,f_s)=\emptyset$.   Set  $d_j:=\deg(f_j)$ and $h_j:=\h(f_j)$ for $1\le j\le s$. Assume that   $d:=d_1\ge d_2\ge \cdots \ge d_{s-1}$ (with no condition with respect to $d_s$), and  $h:=\max\{h_1,\dots, h_{s-1}\}$, and finally set  $r:=\min\{s-1,n\}$.  \\Then there exists $a\in \N\setminus \{0\}$ and $g_1,\dots,g_s\in \Zrng$ such that for $1\le i\le s$ one has
$$a=g_1f_1+\cdots + g_sf_s$$
with $$\deg(g_if_i)\le d_1\cdots d_{r}d_s$$
 and
\begin{align*} \h(a),\h(g_i)+\h(f_i) & \le \Big(\frac{h_s}{d_s} + \sum_{k=1}^{r}\frac{h}{d_k} + (6n+9)\log(n+3)+3n\log\max\{1,s-n\}\Big) d_1\cdots d_{r}d_s\\
&  \le d^{r}h_s + rd^{r-1}d_sh + \big((6n+9)\log(n+3)+3n\log\max\{1,s-n\}\big)d^{r} d_s.
\end{align*}
\end{theorem}

\bigskip
These bounds are again quite tight, as  the following examples show:\\

{\em (1) Brownawell-Masser-Philippon'1986:}  Take
$$f_1=x_1^d,\, f_2=x_1x_n^{d-1}-x_2^d, \dots,\, f_{n-1}=x_{n-2}x_n^{d-1}-x_{n-1}^{d},\, f_{n}=2^h-x_{n-1}x_n^{d-1}.$$
One can show that a B\'ezout identity
$$a=g_1f_1+\cdots +g_{n+1}f_{n+1}$$ specialized at $\bfzeta=(2^{d^{n-2}h}t^{d^{n-1}-1}, \dots, 2^{dh}t^{d^2-1}, 2^ht^{d-1},{1}/{t})$ for any $t\ne 0$ necessarily implies that $a=g_{1}(\bfzeta)2^{
d^{n-1}h}t^{d^n-d}$, which implies that
$$ \deg(g_1f_1) \ge d^n \quad \mbox{and} \quad    \h(a) \ge d^{n-1}h$$
(by specialization at $t=1$ for the second inequality) while the estimates in the previous theorem yield
$$ \max\{\deg(g_if_i)\}\le d^n \quad \mbox{and} \quad    \max \{\h(a),\h(g_i)+\h(f_i) \}  \le  d^{n-1}h + c(n)d^n.$$

{\em (2)} Take the following  modification of the  finite variety $V_2$ introduced in \eqref{eq:V2}:
\begin{equation}\label{eq:NSS2}f_1=x_1-2^h, f_2=x_2-x_1^d, \dots, f_{n}=x_{n}-x_{n-1}^{d}, f_{n+1}=x_n^{d}.\end{equation}
A B\'ezout identity
$$a=g_1f_1+\cdots +g_{n+1}f_{n+1}$$ specialized at  $\bfzeta=(2^h, 2^{dh}, \dots, 2^{d^{n-1}h})$ necessarily implies that $a=g_{n+1}(\bfzeta)2^{
d^nh}$  and therefore $\h(a)\ge d^{n}h$ while the previous theorem yields $$
 \max \{\h(a),\h(g_i)+\h(f_i) \}  \le   d^{n}h + c(n)d^{n+1} .$$
We notice --however-- that this last example doesn't give a simultaneous sharp degree bound: The following  identity  for the case $n=2$ (which readily generalizes to any $n$) $$x_2^d - 2^{d^2h}= \frac{x_2^d - 2^{d^2h}}{x_2-2^{dh}}(x_2-2^{dh})=\frac{x_2^d - 2^{d^2h}}{x_2-2^{dh}}(x_2-x_1^d+x_1^d -2^{dh})$$
implies that $2^{d^2h}= g_1 f_1+ g_2 f_2 + f_3 $ with $\deg (g_1)=2(d-1) $ and $\deg(g_2)=d-1$.\\

The arithmetic Nullstellensatz is often used  in modular methods, when one needs to know if the Nullstellensatz holds over some finite field of characteristic $p$ (see for instance \cite{BBK2009, BGKS2012, CKSZ2014, IKSSS2018}).

\section{An Arithmetic Shape Lemma}

In what follows $I=(f_1,\dots,f_s)\subset \Qrng$ is a {\em radical} zero-dimensional ideal, i.e. $\sqrt I = I$, which  implies that $\dim_\Q\big(\Qrng/I\big)=\deg(V)$  (see for instance \cite[Ch.5, Prop.7]{CLO2015}).\\

The following   Arithmetic {\em Shape Lemma} for the radical zero-dimensional ideal $I$ is an arithmetic version of  what is now also  known as (symbolic) {\em Geometric Resolution}, {\em Kronecker Parameterization} or {\em Rational Univariate Representation}: The classical Shape Lemma already appeared  in the work of  Kronecker \cite{Kro1882}, and it was  reintroduced and adapted to  the context of Computer Algebra around 40 years ago by Chistov and Grigoriev \cite{ChGr1982} and Canny \cite{Can1988}, and again later on, in different forms, by Alonso, Becker et al. \cite[Sec.2.3]{ABRW1996},  Giusti, Heintz et al. \cite{GHMP1995},  and by Fabrice Rouillier \cite{Rou1999} among others. Its meaning is that zero-dimensional radical ideals  are essentially the same as univariate ideals, and working in the quotient $\C[\bfx]/I$ under a suitable isomorphism allows to use the machinery of univariate polynomials. In particular, I will present as an application how this isomorphism allows to obtain improved height estimates  for the representative of a polynomial $p\in \Zrng$ in the quotient algebra $\Qrng/I$.
\\The first part of the following statement is now kind of classical, and the second part presents the height estimates as in \cite[Lem.4.3]{BKM2024} (just slightly more precise due to the application of Theorem~\ref{prop:ABI}). I give a complete proof since it has  its own interest, and also for sake of completeness. As sometimes we know an a priori bound for the degree $\deg(V)$ of the zero-dimensional variety $V$ that  is better than the B\'ezout bound, I leave the height bounds in next statement in terms of this degree. I refer to \cite{GHHMPM1997}, \cite{HMPS2000}, \cite[Ch.13]{Sch2001}, \cite[Prop.4]{SaSc2018}, \cite{HoLe2021} for algorithms, and previous and more general (although less precise) arithmetic estimates.

\begin{theorem} \label{thm:ASL} {\bf (An Arithmetic Shape Lemma)}\\
Let $f_1,\dots,f_s\in \Zrng\setminus \{0\}$ define a {\em radical} zero-dimensional ideal $I\subset \Qrng$, with  $V=V_\C(I)\subset \C^n$ of degree $D:=\deg(V)$.  Set  $d_j:=\deg(f_j)$ and $h_j:=\h(f_j)$ for $1\le j\le s$. Assume that   $d:=d_1\ge d_2\ge \cdots \ge d_{s-1}$ (with no condition with respect to $d_s$), and set $h:=\max\{h_1,\dots, h_{s-1}\}$.  Then there exists an algebra epimorphism
\begin{align*}\varphi: \ \Qrng \quad & \twoheadrightarrow \quad \Q[t]/(\omega_0)\\[1mm]
 x_i\quad & \mapsto \quad {\omega_i(t)}/{{\omega'_0}(t)} \ \mod \ \omega_0 \quad \mbox{for} \ 1\le i\le n\end{align*}
 where  $\omega_0, \omega_1,\dots ,\omega_n\in \Z[t]$, with  $\omega_0$ squarefree, satisfy
 $$\deg(\omega_0)= D, \ \deg(\omega_i)<D \quad \mbox{for} \quad 1\le i\le n$$  and for $0\le i\le n$, $$\h(\omega_i)\le d^{n-1}h_s+ (n-1)d^{n-2}d_sh+2n\log(n+1)d^{n-1}d_s+ 4D\log\big((n+1)D\big).$$
Moreover, the kernel of $\varphi$ satisfies $\ker(\varphi)=I$, and therefore $\varphi$ induces an algebra isomorphism
$$ \overline\varphi: \quad \Qrng/I\quad \displaystyle{\mathop{\longrightarrow}^{\simeq}} \quad  \Q[t]/(\omega_0).$$
\end{theorem}

\begin{proof}

   Let $L(\bfU,\bfx)=U_1x_1+\cdots +U_nx_n\in \Q[\bfU,\bfx]$ be a generic linear form, and consider the polynomial
$$
\Ch_V(t,-\bfU):=\Ch_V(t, -U_1, \ldots, -U_n) = c\prod_{\bfzeta\in V}\big(t-L(\bfU,\bfzeta)\big)\in  \Z[t,\bfU],
$$
where $\Ch_V$ is the (primitive) Chow form of $V$ introduced in Equation~\eqref{eq:chow}.\\
Given  $\bfzeta\in V$,  $\Ch_V\big(L(\bfU,\bfzeta), -\bfU\big)=0$ as a polynomial in $\bfU$ implies by the chain rule that for all $i$  we have
\begin{align*}0&=\partial_{U_i} \big(\Ch_V\big(L(\bfU,\bfzeta), -\bfU\big)\big)(\bfU)\\ & =\partial_t  \big(\Ch_V(t,\bfU)\big) \big(L(\bfU,\bfzeta), -\bfU)\zeta_i- \partial_{U_i}\big(\Ch_V(t,\bfU)\big) (L(\bfU,\bfzeta),-\bfU)\end{align*}
as a polynomial in $\bfU$ too. \\
Therefore, for any $\bfu=(u_1,\dots,u_n)\in \Z^{n}$ such that $\ell(\bfx):=L(\bfu,\bfx)$ satisfies $\ell(\bfzeta)\ne \ell(\bfxi)$ for all $\bfzeta\ne \bfxi \in V$, i.e. $\ell$ is a separating linear form for $V$,
we have that for all $ \bfzeta\in V$, $$ \big(\partial_t \Ch_V(t, -\bfu)\zeta_i- \partial_{U_i} \Ch_V(t,-\bfu)\big)\big(\ell(\bfzeta)\big)=\partial_t \Ch_V(\ell(\bfzeta), -\bfu)\zeta_i- \partial_{U_i} \Ch_V(\ell(\bfzeta),-\bfu)=0.$$ We can then define the univariate integer polynomials $$\omega_0(t):=\Ch_V(t,-\bfu)\quad \mbox{and} \quad \omega_i(t):=\partial_{U_i} \Ch_V(t,-\bfu), \quad \ 1\le i\le n.$$
Since the degree $D$ polynomial $\omega_0 $ has $D$ distinct roots $\ell(\bfzeta),\ \bfzeta\in V$, it is squarefree and its derivative $\omega_0'=\partial_t \Ch_V(t,-\bfu)$ is invertible modulo $\omega_0$.  Moreover,
 for all $\bfzeta\in V$ and $1\le i\le n$,
\begin{equation}\label{eq:zetai}\zeta_i= \dfrac{\partial_{U_i} \Ch_V(\ell(\bfzeta),-\bfu)}{\partial_{t} \Ch_V(\ell(\bfzeta),-\bfu)}=\dfrac{\omega_i(\ell(\bfzeta))}{\omega_0'(\ell(\bfzeta))}.\end{equation}

 This induces a  morphism of algebras $\varphi: \Qrng \to \Q[t]/(\omega_0)$, \ $x_i \mapsto \omega_i(t)/\omega'_0(t) \ \mod \omega_0$, which is well-defined since  $\omega'_0$ is invertible modulo $\omega_0$.
 \\
 Moreover $\varphi$ is an epimorphism because by definition
 $$\varphi(\ell(\bfx))=u_1\varphi(x_1)+\cdots +u_n\varphi(x_n) \equiv u_1\frac{\omega_1(t)}{\omega'_0(t)}+\cdots +u_n\frac{\omega_n(t)}{\omega'_0(t)} \equiv t \quad  \mod  \omega_0$$
since by Identity~\eqref{eq:zetai},  the two polynomials coincide in all the roots $\ell(\bfzeta)$ of the squarefree degree $D$ polynomial $\omega_0\in \Q[t]$:
 $$u_1 \frac{\omega_1(\ell(\bfzeta))}{\omega'_0(\ell(\bfzeta))}+\cdots +u_n\frac{\omega_n(\ell(\bfzeta))}{\omega'_0(\ell(\bfzeta))} =u_1\zeta_1+\cdots +u_n\zeta_n= \ell(\bfzeta).$$

Furthermore,  we show that $\ker(\varphi)=I$:
\begin{align*}g\in \ker(\varphi) &\iff \varphi(g)=0 \iff g\big(\frac{\omega_1(t)}{\omega'_0(t)},\dots,\frac{\omega_n(t)}{\omega'_0(t)}\big) \equiv 0 \quad \mod \omega_0\\
& \iff g\big(\frac{\omega_1(\ell(\bfzeta))}{\omega'_0(\ell(\bfzeta))},\dots,\frac{\omega_n(\ell(\bfzeta))}{\omega'_0(\ell(\bfzeta))}\big) = 0 \quad \mbox{for all } \bfzeta \in V  \quad \mbox{since } \omega_0 \mbox{ is squarefree}\\
&\iff g(\bfzeta)=0  \quad \mbox{for all } \bfzeta \in V  \quad \mbox{by Identity~\eqref{eq:zetai}}\\
&\iff g\in I  \quad \mbox{by the NSS since } I \mbox{ is radical}. \end{align*}
Therefore $\varphi$ induces an isomorphism $$\overline\varphi: \Qrng/I \mathop{\longrightarrow}^{\simeq} \Q[t]/(\omega_0).$$

We now deal with the estimates:

In order to choose $\bfu\in \Z^n$ such that $\ell(\bfx)=u_1x_1+\cdots + u_nx_n\in \Z[\bfx]$  separates the points in $V=\{\bfzeta_1,\dots,\bfzeta_D \}$, we can observe that the non-zero polynomial
$$
\prod_{i<j}(L(\bfU,\bfzeta_i)-L(\bfU,\bfzeta_j))\in \C[\bfU]
$$
has degree $D':=D(D-1)/2< D^2$ and therefore, there exists an element $\bfu=(u_1,\dots,u_n)$ in the integer grid $$\{(k_1,\dots,k_n)\in \Z^{n}\,:\, 0\le k_i\le D'\}$$
where it  does not vanish. In particular, $u_i<D^2$ for  $1\le i\le n$.

By Inequality \eqref{eq:hVChV} and Theorem~\ref{prop:ABI},
\begin{align*}
\h(\Ch_V)&\le \h(V)+3\log(n+1)\deg(V) \\&  \le \ d^{n-1}h_s+ (n-1)d^{n-2}d_sh+2n\log(n+1)d^{n-1}d_s+3D\log(n+1).
\end{align*}
We conclude by observing that the polynomial
$$
\Ch_V(t,-\bfU)=\sum_{\alpha,i:|\alpha|+i=D}a_{\alpha,i}t^i\bfU^\alpha =\sum_{i=0}^D  \big(\sum_{|\alpha|=D-i} a_{\alpha,i}\bfU^\alpha\big)t^i\in \Z[t,\bfU]$$
and therefore \begin{align*}\h(\omega_0)&=\h(\Ch_V(t,-\bfu))\\&\le  \h(\Ch_V)+ D\log(n+1)+ D\log(D^2) \\ & \le d^{n-1}h_s+ (n-1)d^{n-2}d_sh+2n\log(n+1)d^{n-1}d_s+ 4D\log\big((n+1)D\big),\end{align*}
since the number of terms is bounded by $(n+1)^D$.\\
Similarly, for $1\le i\le n$, ${\partial_{U_i} \Ch_V(t,-\bfU)}\in \Z[t,\bfU]$  is homogeneous of degree $D-1$, with height bounded by $ \h(\Ch_V)+\log (D)$, and therefore
\begin{align*}\h(\omega_i)&\le \h(\Ch_V)+ \log(D) + (D-1)\log(n+1) + (D-1)\log(D^2)  \\ &\le d^{n-1}h_s+ (n-1)d^{n-2}d_sh+2n\log(n+1)d^{n-1}d_s+ 4D\log\big((n+1)D\big)\end{align*}
as well.
\end{proof}

Note that in the case of $V_2$ in \eqref{eq:V2}, where there is a single root, the Shape Lemma indicates that $\deg(\omega_0)=1$ while $\deg(\omega_i)=0$ for $1\le i\le n$. In particular, if we take the separating form $\ell=x_1$, we get $\omega_0=t-2^h$ and  $\omega_i=2^{d^{i-1}h}$ for $1\le i\le n$, and a similar bound  holds for any separating linear form $\ell=u_1x_1+\dots+u_nx_n$. This shows that the height bound in the previous theorem might be quite tight. \\

Theorem~\ref{thm:ASL} allows us for instance to get estimates for the heights of the coefficients of a (uniquely defined) representative $\overline p\in \Q[\bfx]/I$ of a polynomial $p\in \Qrng$, that we can call a {\em remainder} of $p$ modulo $I$.  The following follows the developments of \cite[Sec.4.2]{BKM2024}.\\

When $I\subset \Qrng$ is a zero-dimensional radical ideal, then $\dim(\Qrng/I)$ coincides with $D:=\deg(V)\le d^n$, where $d=\max_i\{\deg(f_i)\}$, and $\Qrng/I$ admits a {\em monomial}  basis $\mathcal B$ with $\delta:=\deg({\mathcal B})=\max\{\deg(b) : b\in {\mathcal B}\} < D$: such a basis starts with $b_1:=1$ and then continues for $j>1$ in an inductive process with  $b_j:=x_kb_{i}$, where  $i<j$ and $x_k$ is a variable such that $\{b_1,\dots,b_j\}$ are still linearly independent modulo $I$. Since $\mathcal B$ has $D$ elements, $\delta<D$ (the worst case would be for instance taking $1,x_1,\dots,x_1^{D-1}$).  Moreover, since $b=\bfx^\alpha$ for some $|\alpha|\le \delta$, we have that $\h(b)=0$ for every $b\in {\mathcal B}$.\\

Therefore, for any $p\in \Qrng$, we will get $\overline p=\sum_{b\in {\mathcal B}} c_b b  \,\in \langle {\mathcal B}\rangle_\Q$, with $c_b\in \Q$, as the unique  representative in $\Qrng/I$ of $p\in \Qrng$: this is  what I call the {\em remainder of $p$ modulo $I$ with respect to ${\mathcal B}$}. We observe that there exist of course  $g_1,\dots,g_s\in \Qrng$ such that
\begin{equation}\label{eq:NF}p=g_1f_1+\cdots + g_sf_s+ \overline p.\end{equation}

In order to bound the heights of the numerators and a common denominator for the coefficients of this remainder $\overline p\in \langle {\mathcal B}\rangle_\Q$ we define the following crucial epimorphism ${\mathcal U}$  of $\Q$-vector spaces, which allows us to compute the coefficients in an univariate setting:
$$\begin{array}{cccccccc}{\mathcal U}\ : &  \Qrng & \twoheadrightarrow & \Q[t]/(\omega_0) & \displaystyle{\mathop{\rightarrow}^{\simeq\,}} & \langle 1,t,\dots, t^{D-1}\rangle_\Q\\ &
p& \mapsto &  \omega'_0\,\varphi (p) \ \mod  \omega_0 & \mapsto & (c_0,\dots,c_{D-1})\end{array},$$
where  $\varphi$ is defined in Theorem~\ref{thm:ASL} and $\omega'_0\,\varphi(p)\equiv \sum_{i=0}^{D-1}c_it^i \ \mod \omega_0$.
Observe that since $\omega'_0 $ is invertible modulo $\omega_0$,  $\mathcal U$ is a $\Q$-epimorphism with $\ker({\mathcal U})=\ker (\varphi)=I$.

\begin{lemma} {\em (\cite[Lem.4.4]{BKM2024})}\\Let $f_1,\dots,f_s\in \Zrng\setminus \{0\}$ define a  radical zero-dimensional ideal $I\subset \Qrng$ with $V=V_\C(I)$, and let $\bfx^\alpha$ be a monomial. Set  $d:=\max\{\deg(f_j): 1\le j\le s\}$ and $h:=\max\{\h(f_j): 1\le j\le s\}$.   Then
$$\h({\mathcal U}(\bfx^\alpha))\le  nd^{n-1}|\alpha|h+ 4n\log((n+2)d)d^n|\alpha|.$$
\end{lemma}

\begin{proof} Let $J:=\big(f_1,\dots,f_s,x_{n+1}-\bfx^\alpha\big)\subset \Q[\bfx,x_{n+1}]$, which is a radical  zero-dimensional ideal of degree $\deg(V)\le d^n$, since $V_\C(J)=\{(\bfzeta, \bfzeta^\alpha):\, \bfzeta\in V\}$ as already seen before. Note that if $\ell(\bfx)=u_1x_1+\cdots + u_nx_n\in \Zrng$ is a separating linear form for $V$, it is still a separating linear form for $V_\C(J)$. Therefore in the construction of $\varphi_J:\Q[\bfx,x_{n+1}]\to \Q[t]/(\omega^J_0)$ for the ideal $J$ in Theorem~\ref{thm:ASL}, we can choose this linear form $\ell(\bfx)$, which also applies for $\varphi_I:\Qrng\to \Q[t]/(\omega^I_0)$.\\
 We observe that since $u_{n+1}=0$ in the choice of the linear form, we have that $\omega_i^J=\omega^I_i$ for $0\le i\le n$, which implies that $\varphi_J(x_i)=\varphi_I(x_i)$ for $1\le i\le n$.  Let  $\omega_{n+1}:=\omega_{n+1}^J(t)\in \Z[t]$ be such that
$\varphi_J(x_{n+1})=(\omega'_0)^{-1}\omega_{n+1}$.\\
Then, since  $x_{n+1}\equiv \bfx^\alpha\ \mod J$, we have
$${\mathcal U}(\bfx^\alpha)=\omega'_0 \varphi_I(\bfx^\alpha) =\omega'_0 \varphi_J(\bfx^\alpha) =\omega'_0\varphi_J(x_{n+1})=\omega_{n+1}.$$
Finally, applying Theorem~\ref{thm:ASL} to $J\subset \Q[\bfx,x_{n+1}]$ with $\deg(J)=\deg(I) \le d^n$ and $f_{s+1}=x_{n+1}-\bfx^\alpha$ of degree $|\alpha|$ and height $0$, we get
\begin{align*}\h({\mathcal U}(\bfx^\alpha))=\h(\omega_{n+1})&\le   nd^{n-1}|\alpha|h+ 2(n+1)\log(n+2)d^n|\alpha|+ 4nd^n\log((n+1)d)\\ & \le
nd^{n-1}|\alpha|h+ 4n\log((n+2)d)d^n|\alpha|.\end{align*}
\end{proof}

\begin{corollary} {\em (\cite[Cor.4.5]{BKM2024})}\\Let $f_1,\dots,f_s\in \Zrng \setminus \{0\}$ define a  radical zero-dimensional ideal $I\subset \Qrng$ and let $p\in \Zrng\setminus \{0\}$. Set  $d:=\max\{\deg(f_j): 1\le j\le s\}$, $d_p:=\deg(p)$, and $h:=\max\{\h(f_j): 1\le j\le s\}$, $h_p:=\h(p)$.     Then
$$\h({\mathcal U}(p))\le   h_p+  nd^{n-1}d_ph+ 5n\log((n+2)d)d^nd_p.$$
\end{corollary}

\begin{proof} This is as the proof of Corollary~\ref{cor:hproot}: Write $p=\sum_\alpha p_\alpha \bfx^\alpha $ with $p_\alpha\in \Z$. Then $${\mathcal U}(p)=\sum_\alpha p_\alpha {\mathcal U}(\bfx^\alpha)\quad \mbox{
with } \quad \h({\mathcal U}(\bfx^\alpha))\le  nd^{n-1}d_ph+ 4n\log((n+2)d)d^nd_p.$$
Since $p$ has at most $(d_p+1)^n$ monomials, we have
\begin{align*}\h({\mathcal U}(p))&\le h_p+ nd^{n-1}d_ph+ 4n\log((n+2)d)d^nd_p+ n\log(d_p+1) \\ & \le h_p+  nd^{n-1}d_ph+ 5n\log((n+2)d)d^nd_p.\end{align*}

\end{proof}

We are now ready to prove a height bound as presented in  \cite[Prop.4.6]{BKM2024} for the numerator and a common denominator of the remainder $\overline p\in \langle \mathcal{B}\rangle_\Q$ of an integer polynomial $p\in \Zrng$ modulo $I$ (except that here we do not need to assume that $\delta=\deg(B)\le d$ and $\deg(p)\ge d$).
\begin{proposition} {\bf (Height of the remainder modulo $I$)}\\
Let $f_1,\dots,f_s\in \Zrng\setminus \{0\}$ define a  radical zero-dimensional ideal $I\subset \Qrng$,  ${\mathcal B}$ be a monomial basis of $\Q[\bfx]/I$ with $\delta:=\deg({\mathcal B})$ and   let $p\in \Zrng\setminus \{0\}$.  Set  $d:=\max\{\deg(f_j): 1\le j\le s\}$, $d_p:=\deg(p)$, and $h:=\max\{\h(f_j): 1\le j\le s\}$, $h_p:=\h(p)$.   \\Then there exist   $a\in \N$ and $N(p)\in \langle  \mathcal{B}\rangle_\Z$ such that $\overline p= N(p)/a$ with
\begin{align*}&\h(a)\le  nd^{2n-1}\delta h+ 5n\log((n+2)d)d^{2n}\delta,\\
& \h({N(p)})\le h_p+ nd^{n-1}(d_p+d^n\delta)h + 5n\log((n+2)d)d^n(d_p+d^n\delta).\end{align*}
\end{proposition}
\begin{proof}
If $\overline p=\sum_{b\in {\mathcal B}}c_b \, b$ with $c_b\in \Q$, we have that  $${\mathcal U}(p)={\mathcal U}(\overline p)=\sum_{b\in {\mathcal B}}c_b \, {\mathcal U}(b)$$ since $p\equiv \overline p \ \mod I$.
Therefore $(c_b: b\in {\mathcal B})$ is the (unique) solution $y$ of a square  linear system of equations  $$A\,y=c$$ of size $D$,  where the matrix $A\in \Z^{D\times D}$  to invert is composed by the coefficients of ${\mathcal U}(b)$, $b\in {\mathcal B}$, and the constant vector $c$ is composed by the integer coefficients of ${\mathcal U}(p)$ in the basis $(1,t,\dots, t^{D-1})$.
Solving this system of equations by Cramer's rule gives the common denominator $a=\det(A)$ and the coordinates of the numerator  $N(p)\in \Zrng$ where each column of $A$ is replaced by the constant vector. \\
By Hadamard's bound, the determinant of a  matrix with columns $v_1,\dots, v_D\in \Z^D$ where all but one column has height bounded by $h_A$ and one column has height bounded by $h_c$ satisfies
$$\h(\det(v_1,\dots, v_D))\le h_c+(D-1)h_A + \frac{D\log(D)}{2}.$$
Since $$h_A\le  nd^{n-1}\delta h+ 4n\log((n+2)d)d^n\delta \quad \mbox{and} \quad  h_c\le h_p+  nd^{n-1}d_ph+ 5n\log((n+2)d)d^nd_p$$ we  conclude that
\begin{align*} \h(a) & \le D\big( nd^{n-1}\delta h+ 4n\log((n+2)d)d^n\delta\big) +\frac{D\log(D)}{2}\\
& \le nd^{2n-1}\delta h+ 5n\log((n+2)d)d^{2n}\delta  \end{align*}
and
\begin{align*}\h(N(p))&\le h_p+  nd^{n-1}d_ph+ 5n\log((n+2)d)d^nd_p\\ & \qquad \ +(D-1)\big( nd^{n-1}\delta h+ 4n\log((n+2)d)d^n\delta\big) +\frac{D\log(D)}{2}\\
& \le h_p+ nd^{n-1}(d_p+d^n\delta)h + 5n\log((n+2)d)d^n(d_p+d^n\delta).\end{align*}
\end{proof}
We observe that the same bound for $\h(\overline p)$ holds for the maximum  of the logarithms of the absolute values of  the coefficients of the polynomial $\overline p$ when $p\in \C[\bfx]$. This is because the common denominator $a$ satisfies $|a|\ge 1$:

\begin{corollary} Let $f_1,\dots,f_s\in \Zrng\setminus \{0\}$ define a  radical zero-dimensional ideal $I\subset \Crng$,  ${\mathcal B}$ be a monomial basis of $\C[\bfx]/I$ with $\delta:=\deg({\mathcal B})$ and   let $p\in \Crng\setminus \{0\}$.  Set  $d:=\max\{\deg(f_j): 1\le j\le s\}$, $d_p:=\deg(p)$, and $h:=\max\{\h(f_j): 1\le j\le s\}$, $h_p$ for  the maximum  of the logarithms of the absolute values of  the coefficients of $p$. \\ Then $\overline p=\sum_{b\in B} c_b \,b\in \langle\mathcal{B}\rangle_\C$ satisfies
$$\max\{\log |c_b| :\, b\in B\} \le  h_p+ nd^{n-1}(d_p+d^n\delta)h + 5n\log((n+2)d)d^n(d_p+d^n\delta).$$
\end{corollary}

This estimate is essentially better than the bound one would obtain by applying directly a linear algebra argument as in the previous proof, even in the case when the polynomials $f_1,\dots,f_s$ are a  degree-preserving Gr\"obner basis of $I$ (for any $p\in I$, one has that $p=g_1f_1+\cdots + g_sf_s$ with
$\deg(g_j)\le \deg(p)-\deg(f_j)$, for $1\le j\le s$),   since for $p$ of large degree $d_p$, we would have to solve a linear system of order $d_p^n$.\\

Let me close this section by noticing that this estimate might be again (quite)  tight in terms of the height, as suggested by  Example~\eqref{eq:NSS2}: For $$f_1:=x_1-2^h,f_2:=x_2-x_1^d,\dots, f_n=x_n-x_{n-1}^d \quad \mbox{and} \quad p=x_n^d$$ we have $\mathcal{B}=\{1\}$ and $\overline p=2^{d^nh}$.

\section{An Arithmetic Perron's Theorem}

In \cite{Jel2005} the proof of the effective Nullstellensatz relies on bounding the degrees of a coefficient in a given algebraic equation for polynomials related to the polynomials $f_1,\dots,f_s$, which define the empty variety.
I present here a sketch of an algebraic version of his proof for the case $s=n+1$,  which is the typical case (the same  holds for $s< n+1$ and for $s>n+1$ one usually performs some linear combinations of $f_1,\dots,f_s$  to restrict to $n+1$ polynomials). \\
Assume that $f_1,\dots, f_{n+1}\in \Crng$ satisfy that $V_\C(f_1,\dots,f_{s})=\emptyset$, so that by the NSS there exist (unknown) $g_1,\dots, g_{n+1}$ with
$$1=g_1f_1+\dots + g_{n+1}f_{n+1},$$ and consider the following algebra morphism
$$\Phi: \ \C[\bfx,z_1,\dots,z_{n+1}] \  \rightarrow \ \C[\bfx,z], \
 x_i\  \mapsto \ {x_i}, \
  z_j\ \mapsto \\ {zf_j(\bfx)} \  \mbox{ for} \ 1\le i\le n, 1\le j\le n+1.
  $$
The map $\Phi$ turns out to be an epimorphism because $\Phi(g_1(\bfx)z_1+\cdots + g_{n+1}(\bfx)z_{n+1}) = z$. Therefore,   $\Phi$ induces an algebra isomorphism $$\overline\Phi: \ A:=\C[\bfx,z_1,\dots,z_{n+1}]/\ker(\Phi)\  \simeq \  \C[\bfx,z],$$ between the finitely generated algebra $A$ over $\C$ and the polynomial ring $\C[\bfx,z]$,
and in particular  the (Krull) dimension  of $A$ equals $n+1$.

\smallskip

By Noether's normalization lemma (see for instance \cite[Ch.5, Ex.16]{AtMa1969}) there exist  $n+1$ linear combinations of the variables $\bfx,z_1,\dots,z_{n+1}$ in $A$, that can moreover be taken of the form    $z_1+\ell_1(\bfx),\dots,z_{n+1}+\ell_{n+1}(\bfx)$ (where $\ell_1,\dots, \ell_{n+1}$ are linear forms in the variables $\bfx$)  that are algebraically independent  over $\C$  and such that $A$ is integral over $\C[z_1+\ell_1(\bfx),\dots, z_{n+1}+\ell_{n+1}(\bfx)]$.

\smallskip

In particular,  $g_1(\bfx)z_1+\dots+g_{n+1}(\bfx)z_{n+1}$ is integral over $\C[z_1+\ell_1(\bfx),\dots, z_{n+1}+\ell_{n+1}(\bfx)]$: Setting $\bfy:=(y_1,\dots,y_{n+1})$, there exists a polynomial of some degree $N$,
$$P(\bfy, t)=t^N +\sum_{j=1}^{N}\sum_{\alpha} a_{\alpha,j}\bfy^\alpha t^{N-j}\quad \in \, \C[\bfy,t]\setminus \{0\}, $$
monic in $t$,
 such that
$$P\big(z_1+\ell_1(\bfx), \dots, z_{n+1}+\ell_{n+1}(\bfx), g_1(\bfx)z_1+\dots+g_{n+1}(\bfx)z_{n+1}\big) \,\equiv \, 0 \quad \mod \ker(\Phi).$$
Therefore
\begin{align*} 0& =\Phi\big(P\big(z_1+\ell_1(\bfx), \dots, z_{n+1}+\ell_{n+1}(\bfx), g_1(\bfx)z_1+\dots+g_{n+1}(\bfx)z_{n+1}\big)\big)\\
& = P \big( \Phi\big(z_1+\ell_1(\bfx)\big),\dots,\Phi\big(z_{n+1}+\ell_{n+1}(\bfx)\big), \Phi\big( g_1(\bfx)z_1+\dots+g_{n+1}(\bfx)z_{n+1}\big) \big)\\
& = P\big(zf_1(\bfx)+ \ell_1(\bfx),\dots, zf_{n+1}(\bfx)+\ell_{n+1}(\bfx),z \big) \\
& = z^N +\sum_{j=1}^{N}\sum_\alpha a_{\alpha,j}\big(zf_1(\bfx)+ \ell_1(\bfx)\big)^{\alpha_1}\cdots  \big(zf_{n+1}(\bfx)+ \ell_{n+1}(\bfx)\big)^{\alpha_{n+1}} \,z^{N-j}.\end{align*}
By inspection, the coefficient of $z^N$ when expanding this expression
gives an effective B\'ezout identity for $f_1,\dots,f_s$.\\

Finding the degrees in this integral dependence equation  for  $z$ is strongly related to the problem of finding an  algebraic dependence equation for polynomials that are not algebraically independent.
This is what is sometimes called {\em Perron's theorem}, since the very popular book by Perron \cite{Per1927} gives in Satz 57 a combinatorial proof that if
$f_1,\dots,f_{n+1}\in \C[\bfx]$ are $n+1$ polynomials in $n$ variables of degrees $d_j:=\deg(f_j)$,  then not only they are algebraically dependent, but  there exists a polynomial $$P=\sum_\alpha c_\alpha y_1^{\alpha_1}\cdots y_{n+1}^{\alpha_n}\quad \in \C[y_1,\dots,y_{n+1}]\setminus \{0\}$$
which satisfies that
\begin{itemize}\item $P(f_1,\dots,f_{n+1})=0$,
\item $\alpha_1d_1+\cdots +\alpha_{n+1}d_{n+1}\le d_1\cdots d_{n+1}$ \ for all $\alpha$.
\end{itemize}

\smallskip

In \cite[Cor.3.23]{DKS2013} we got the following arithmetic version of Perron's theorem as a consequence of a more general version over a variety. A parameterized version of it was a crucial tool to obtain the bounds in the Arithmetic Nullstellensatz presented in Theorem~\ref{thm:ANSS}.

\begin{theorem} {\bf (An arithmetic Perron's theorem)}\\
Let $f_1,\dots,f_{n+1}\in \Zrng\setminus \{0\}$, and set   $d_j:=\deg(f_j)$ and $h_j:=\h(f_j)$ for $1\le j\le n+1$. Then there exists a polynomial $$P=\sum_\alpha c_\alpha y_1^{\alpha_1}\dots y_{n+1}^{\alpha_{n+1}}\quad \in \, \Z[y_1,\dots,y_{n+1}]\setminus \{0\}$$  which satisfies that  \begin{itemize} \item $P(f_1,\dots,f_{n+1})=0$,
 \item $\alpha_1d_1+\cdots +\alpha_{n+1}d_{n+1}\le d_1\cdots d_{n+1}$ \ for all $\alpha$,
\item $\h(c_\alpha)+ \displaystyle\sum_{i=1}^{n+1}\alpha_i h_i\le \Big(\displaystyle \sum_{i=1}^{n+1}\dfrac{h_i}{d_i} + (n+2)\log(2n+8)\Big) d_1\cdots d_{n+1}$ \  for all $\alpha$ s.t. $c_\alpha\ne 0$.
\end{itemize}
\end{theorem}
I end this text by showing  that these bounds appear again to be essentially tight,  as can be shown by the same Example~\eqref{eq:NSS2} given for the arithmetic Nullstellensatz:
\begin{example}\label{ex:perron} Let $$f_1=x_1-2^h, f_2=x_2-x_1^d, \dots, f_{n}=x_{n}-x_{n-1}^{d}, f_{n+1}=x_n^{d}\ \in \Z[\bfx]$$ and let $P(y_1,\dots,y_{n+1})\in \Z[y_1,\dots,y_{n+1}]$ be a (non-zero) algebraic dependence equation for them. Then $$\deg_{y_{1}}(P)\ge d^n\quad \mbox{and}\quad  \h(P)\ge d^nh.$$
\end{example}
\begin{proof} Here is an elementary proof, although surely not the shortest one! \\We first show that  $f_1$ is integral over $\Z[f_2,\dots,f_{n+1}]$ by induction on $n$:\\
For $n=1$, we have that $f_1=x-2^h$ and $f_2=x^d$: since $2^h\in \Z$ is integral over $\Z[f_2]$, it is enough to verify that $x$ is integral over $\Z[f_2]$, which is clear since $x^d-f_2=0$. \\
To prove the statement for $n$, we observe that by the inductive hypothesis, $f_1$ is integral over $\Z[f_2,\dots,f_{n-1},x_{n-1}^d]$ and also that
$x_{n-1}^d $ is integral over $\Z[f_n,f_{n+1}]$ since
$$x_n^{d}=\big((x_n -x_{n-1}^d)+ x_{n-1}^d\big)= \sum_{k=0}^d\binom{d}{k} (x_n-x_{n-1}^d)^{d-k}(x_{n-1}^{d})^k  $$
implies that $x_{n-1}^d$ is a root of the monic polynomial
$$y^d+\sum_{k=0}^{d-1}\binom{d}{k}f_{n}^{d-k}y^{k} - f_{n+1}\ \in \Z[f_{n},f_{n+1}][y]. $$
Therefore, $f_1$  is integral over  $\Z[f_2,\dots,f_{n+1}]$.\\
 Let $Q\in \Q(y_2,\dots,y_{n+1})[y_1]$ be such that $Q(f_2,\dots,f_{n+1})[y_1]$ is  the minimal (monic) polynomial  of $f_1$ over $\Q(f_2,\dots,f_{n+1})$. The fact that $f_1$ is integral over $\Z[f_2,\dots,f_{n+1}]$ then implies  that  $Q\in \Z[y_1,\dots,y_{n+1}]$.
Write
$$Q=q(y_1,y_{n+1}) + y_2\,q_2(\bfy)+\cdots + y_n\,q_n(\bfy),$$
where $q\ne 0$ since $Q$ is monic in $y_1$.
We have $$Q(f_1,\dots,f_{n+1})=q(x_1-2^h,x_n^d)+ (x_2-x_1^d)q_2(\star)+\dots+ (x_n-x_{n-1}^d)q_n(\star)=0$$ and by specializing this expression at
$(x_1,x_1^{d},\dots,x_1^{d^{n-1}})$ we deduce that
$$q(x_1-2^h,x_1^{d^n})=0.$$
This implies that $y_{n+1}-(y_1+2^h)^{d^n} \mid q(y_1,y_{n+1} ) $ in $\Z[y_1,y_{n+1}]$, since by polynomial division, the identity
$$q(y_1,y_{n+1}) =q_1(y_1,y_{n+1})\big(y_{n+1}-(y_1+2^h)^{d^n}\big)+ r(y_1)$$ specialized at $y_1=x_1-2^h, y_{n+1}= (y_1+2^h)^{d^n}= x_1^{d^n}$ yields $r(x_1-2^h)=0$ as a polynomial, i.e. $r=0$.\\
Now, from  $q\ne 0 $ we deduce that $\deg_{y_1}(q) \ge d^n$ and $\h(q)\ge d^nh$, by looking at its constant term. Therefore the same happens with $Q$ since $q$ cannot be canceled by the other terms in $Q$, which are all divisible by some $y_j$, $2\le j\le  n$. \\ Finally, any other polynomial dependence equation  $P\in \Z[y_1,\dots,y_{n+1}]$ for $f_1,\dots, f_{n+1}$ has to involve the variable $y_1$ since  $f_2,\dots, f_{n+1}$ are algebraically independent over $\Q$ (any algebraic dependence equation for them gives an algebraic dependence equation for $x_1,\dots,x_n$) and therefore is divisible in $\Z[y_1,\dots,y_{n+1}]$ by $Q$. This concludes the proof.
\end{proof}

\end{document}